\documentstyle[12pt]{article}

\expandafter\chardef\csname pre amssym.def at\endcsname=\the\catcode`\@
\catcode`\@=11

\def\undefine#1{\let#1\undefined}
\def\newsymbol#1#2#3#4#5{\let\next@\relax
 \ifnum#2=\@ne\let\next@\msafam@\else
 \ifnum#2=\tw@\let\next@\msbfam@\fi\fi
 \mathchardef#1="#3\next@#4#5}
\def\mathhexbox@#1#2#3{\relax
 \ifmmode\mathpalette{}{\m@th\mathchar"#1#2#3}%
 \else\leavevmode\hbox{$\m@th\mathchar"#1#2#3$}\fi}
\def\hexnumber@#1{\ifcase#1 0\or 1\or 2\or 3\or 4\or 5\or 6\or 7\or 8\or
 9\or A\or B\or C\or D\or E\or F\fi}

\font\tenmsa=msam10 scaled \magstep1
\font\sevenmsa=msam7 scaled \magstep1
\font\fivemsa=msam5 scaled \magstep1
\newfam\msafam
\textfont\msafam=\tenmsa
\scriptfont\msafam=\sevenmsa
\scriptscriptfont\msafam=\fivemsa
\edef\msafam@{\hexnumber@\msafam}
\mathchardef\dabar@"0\msafam@39
\def\dashrightarrow{\mathrel{\dabar@\dabar@\mathchar"0\msafam@4B}}
\def\dashleftarrow{\mathrel{\mathchar"0\msafam@4C\dabar@\dabar@}}

\def\ulcorner{\delimiter"4\msafam@70\msafam@70 }
\def\urcorner{\delimiter"5\msafam@71\msafam@71 }
\def\llcorner{\delimiter"4\msafam@78\msafam@78 }
\def\lrcorner{\delimiter"5\msafam@79\msafam@79 }
\def\yen{{\mathhexbox@\msafam@55 }}
\def\checkmark{{\mathhexbox@\msafam@58 }}
\def\circledR{{\mathhexbox@\msafam@72 }}
\def\maltese{{\mathhexbox@\msafam@7A }}

\font\tenmsb=msbm10 scaled \magstep1
\font\sevenmsb=msbm7 scaled \magstep1
\font\fivemsb=msbm5 scaled \magstep1
\newfam\msbfam
\textfont\msbfam=\tenmsb
\scriptfont\msbfam=\sevenmsb
\scriptscriptfont\msbfam=\fivemsb
\edef\msbfam@{\hexnumber@\msbfam}
\def\Bbb#1{\fam\msbfam\relax#1}
\def\widehat#1{\setboxz@h{$\m@th#1$}%
 \ifdim\wdz@>\tw@ em\mathaccent"0\msbfam@5B{#1}%
 \else\mathaccent"0362{#1}\fi}
\def\widetilde#1{\setboxz@h{$\m@th#1$}%
 \ifdim\wdz@>\tw@ em\mathaccent"0\msbfam@5D{#1}%
 \else\mathaccent"0365{#1}\fi}
\font\teneufm=eufm10 scaled \magstep1
\font\seveneufm=eufm7 scaled \magstep1
\font\fiveeufm=eufm5 scaled \magstep1
\newfam\eufmfam
\textfont\eufmfam=\teneufm
\scriptfont\eufmfam=\seveneufm
\scriptscriptfont\eufmfam=\fiveeufm
\def\frak#1{{\fam\eufmfam\relax#1}}

\catcode`\@=\csname pre amssym.def at\endcsname

\expandafter\ifx\csname pre amssym.tex at\endcsname\relax \else \endinput\fi
\expandafter\chardef\csname pre amssym.tex at\endcsname=\the\catcode`\@
\catcode`\@=11
\newsymbol\boxdot 1200
\newsymbol\boxplus 1201
\newsymbol\boxtimes 1202
\newsymbol\square 1003
\newsymbol\blacksquare 1004
\newsymbol\centerdot 1205
\newsymbol\lozenge 1006
\newsymbol\blacklozenge 1007
\newsymbol\circlearrowright 1308
\newsymbol\circlearrowleft 1309
\undefine\rightleftharpoons
\newsymbol\rightleftharpoons 130A
\newsymbol\leftrightharpoons 130B
\newsymbol\boxminus 120C
\newsymbol\Vdash 130D
\newsymbol\Vvdash 130E
\newsymbol\vDash 130F
\newsymbol\twoheadrightarrow 1310
\newsymbol\twoheadleftarrow 1311
\newsymbol\leftleftarrows 1312
\newsymbol\rightrightarrows 1313
\newsymbol\upuparrows 1314
\newsymbol\downdownarrows 1315
\newsymbol\upharpoonright 1316
 \let\restriction\upharpoonright
\newsymbol\downharpoonright 1317
\newsymbol\upharpoonleft 1318
\newsymbol\downharpoonleft 1319
\newsymbol\rightarrowtail 131A
\newsymbol\leftarrowtail 131B
\newsymbol\leftrightarrows 131C
\newsymbol\rightleftarrows 131D
\newsymbol\Lsh 131E
\newsymbol\Rsh 131F
\newsymbol\rightsquigarrow 1320
\newsymbol\leftrightsquigarrow 1321
\newsymbol\looparrowleft 1322
\newsymbol\looparrowright 1323
\newsymbol\circeq 1324
\newsymbol\succsim 1325
\newsymbol\gtrsim 1326
\newsymbol\gtrapprox 1327
\newsymbol\multimap 1328
\newsymbol\therefore 1329
\newsymbol\because 132A
\newsymbol\doteqdot 132B
 
\newsymbol\triangleq 132C
\newsymbol\precsim 132D
\newsymbol\lesssim 132E
\newsymbol\lessapprox 132F
\newsymbol\eqslantless 1330
\newsymbol\eqslantgtr 1331
\newsymbol\curlyeqprec 1332
\newsymbol\curlyeqsucc 1333
\newsymbol\preccurlyeq 1334
\newsymbol\leqq 1335
\newsymbol\leqslant 1336
\newsymbol\lessgtr 1337
\newsymbol\backprime 1038
\newsymbol\risingdotseq 133A
\newsymbol\fallingdotseq 133B
\newsymbol\succcurlyeq 133C
\newsymbol\geqq 133D
\newsymbol\geqslant 133E
\newsymbol\gtrless 133F
\newsymbol\sqsubset 1340
\newsymbol\sqsupset 1341
\newsymbol\vartriangleright 1342
\newsymbol\vartriangleleft 1343
\newsymbol\trianglerighteq 1344
\newsymbol\trianglelefteq 1345
\newsymbol\bigstar 1046
\newsymbol\between 1347
\newsymbol\blacktriangledown 1048
\newsymbol\blacktriangleright 1349
\newsymbol\blacktriangleleft 134A
\newsymbol\vartriangle 134D
\newsymbol\blacktriangle 104E
\newsymbol\triangledown 104F
\newsymbol\eqcirc 1350
\newsymbol\lesseqgtr 1351
\newsymbol\gtreqless 1352
\newsymbol\lesseqqgtr 1353
\newsymbol\gtreqqless 1354
\newsymbol\Rrightarrow 1356
\newsymbol\Lleftarrow 1357
\newsymbol\veebar 1259
\newsymbol\barwedge 125A
\newsymbol\doublebarwedge 125B
\undefine\angle
\newsymbol\angle 105C
\newsymbol\measuredangle 105D
\newsymbol\sphericalangle 105E
\newsymbol\varpropto 135F
\newsymbol\smallsmile 1360
\newsymbol\smallfrown 1361
\newsymbol\Subset 1362
\newsymbol\Supset 1363
\newsymbol\Cup 1264
 
\newsymbol\Cap 1265
 
\newsymbol\curlywedge 1266
\newsymbol\curlyvee 1267
\newsymbol\leftthreetimes 1268
\newsymbol\rightthreetimes 1269
\newsymbol\subseteqq 136A
\newsymbol\supseteqq 136B
\newsymbol\bumpeq 136C
\newsymbol\Bumpeq 136D
\newsymbol\lll 136E
 
\newsymbol\ggg 136F
 
\newsymbol\circledS 1073
\newsymbol\pitchfork 1374
\newsymbol\dotplus 1275
\newsymbol\backsim 1376
\newsymbol\backsimeq 1377
\newsymbol\complement 107B
\newsymbol\intercal 127C
\newsymbol\circledcirc 127D
\newsymbol\circledast 127E
\newsymbol\circleddash 127F
\newsymbol\lvertneqq 2300
\newsymbol\gvertneqq 2301
\newsymbol\nleq 2302
\newsymbol\ngeq 2303
\newsymbol\nless 2304
\newsymbol\ngtr 2305
\newsymbol\nprec 2306
\newsymbol\nsucc 2307
\newsymbol\lneqq 2308
\newsymbol\gneqq 2309
\newsymbol\nleqslant 230A
\newsymbol\ngeqslant 230B
\newsymbol\lneq 230C
\newsymbol\gneq 230D
\newsymbol\npreceq 230E
\newsymbol\nsucceq 230F
\newsymbol\precnsim 2310
\newsymbol\succnsim 2311
\newsymbol\lnsim 2312
\newsymbol\gnsim 2313
\newsymbol\nleqq 2314
\newsymbol\ngeqq 2315
\newsymbol\precneqq 2316
\newsymbol\succneqq 2317
\newsymbol\precnapprox 2318
\newsymbol\succnapprox 2319
\newsymbol\lnapprox 231A
\newsymbol\gnapprox 231B
\newsymbol\nsim 231C
\newsymbol\ncong 231D
\newsymbol\diagup 231E
\newsymbol\diagdown 231F
\newsymbol\varsubsetneq 2320
\newsymbol\varsupsetneq 2321
\newsymbol\nsubseteqq 2322
\newsymbol\nsupseteqq 2323
\newsymbol\subsetneqq 2324
\newsymbol\supsetneqq 2325
\newsymbol\varsubsetneqq 2326
\newsymbol\varsupsetneqq 2327
\newsymbol\subsetneq 2328
\newsymbol\supsetneq 2329
\newsymbol\nsubseteq 232A
\newsymbol\nsupseteq 232B
\newsymbol\nparallel 232C
\newsymbol\nmid 232D
\newsymbol\nshortmid 232E
\newsymbol\nshortparallel 232F
\newsymbol\nvdash 2330
\newsymbol\nVdash 2331
\newsymbol\nvDash 2332
\newsymbol\nVDash 2333
\newsymbol\ntrianglerighteq 2334
\newsymbol\ntrianglelefteq 2335
\newsymbol\ntriangleleft 2336
\newsymbol\ntriangleright 2337
\newsymbol\nleftarrow 2338
\newsymbol\nrightarrow 2339
\newsymbol\nLeftarrow 233A
\newsymbol\nRightarrow 233B
\newsymbol\nLeftrightarrow 233C
\newsymbol\nleftrightarrow 233D
\newsymbol\divideontimes 223E
\newsymbol\varnothing 203F
\newsymbol\nexists 2040
\newsymbol\Finv 2060
\newsymbol\Game 2061
\newsymbol\mho 2066
\newsymbol\eth 2067
\newsymbol\eqsim 2368
\newsymbol\beth 2069
\newsymbol\gimel 206A
\newsymbol\daleth 206B
\newsymbol\lessdot 236C
\newsymbol\gtrdot 236D
\newsymbol\ltimes 226E
\newsymbol\rtimes 226F
\newsymbol\shortmid 2370
\newsymbol\shortparallel 2371
\newsymbol\smallsetminus 2272
\newsymbol\thicksim 2373
\newsymbol\thickapprox 2374
\newsymbol\approxeq 2375
\newsymbol\succapprox 2376
\newsymbol\precapprox 2377
\newsymbol\curvearrowleft 2378
\newsymbol\curvearrowright 2379
\newsymbol\digamma 207A
\newsymbol\varkappa 207B
\newsymbol\Bbbk 207C
\newsymbol\hslash 207D
\undefine\hbar
\newsymbol\hbar 207E
\newsymbol\backepsilon 237F
\catcode`\@=\csname pre amssym.tex at\endcsname

\newcommand{\ga}{\alpha}     
\newcommand{\gb}{\beta}      
\renewcommand{\gg}{\gamma}   
\newcommand{\gd}{\delta}

\newcommand{\gth}{\theta}    
      
\newcommand{\gk}{\kappa}  
\newcommand{\gl}{\lambda}    
\newcommand{\gm}{\mu}        
\newcommand{\gn}{\nu}

\newcommand{\gr}{\rho}       
\newcommand{\gs}{\sigma}     
\newcommand{\gt}{\tau}

\newcommand{\go}{\omega}

\newcommand{\gD}{\Delta}

\newcommand{\ha}{\aleph}

\newcommand{\setof}[2]{{\{\; #1 \; \vert \; #2 \; \} } }
\newcommand{\seq}[1]{{\langle #1 \rangle} }
\newcommand{\card}[1]{{\vert #1 \vert} }

\newcommand{\forces}{\Vdash}

\renewcommand{\models}{\vDash}
\newcommand{\powerset}{{\cal P}}
\newcommand{\dom}{{\rm dom}}
\newcommand{\rge}{{\rm rge}}

\newcommand{\cf}{{\rm cf}}
\newcommand{\lh}{{\rm lh}}

\newcommand{\lra}{\longrightarrow}

\newtheorem{definition}{Definition}
\newtheorem{theorem}{Theorem}
\newenvironment{proof}{\noindent{\bf Proof:}}{\nopagebreak\mbox{}\newline
 \makebox[\textwidth]{\hfill$\blacklozenge$}\par\bigskip}

\newcommand{\implies}{\Longrightarrow}

\catcode`\@=11
\def\@begintheorem#1#2{\rm \trivlist \item[\hskip \labelsep{\bf #1\ #2:}]}
\def\@opargbegintheorem#1#2#3{\rm \trivlist
      \item[\hskip \labelsep{\bf #1\ #2\ (#3):}]}
\catcode`\@=12

\newtheorem{lemma}{Lemma}

\newtheorem{claim}{Claim}

\newtheorem{fact}{Fact}

\newcommand{\FP}{{\Bbb P}}
\newcommand{\FQ}{{\Bbb Q}}
\newcommand{\FR}{{\Bbb R}}

\newcommand{\FD}{{\Bbb D}}

\title{Cardinal invariants above the continuum}

\author{James Cummings\\
        Hebrew University of Jerusalem\\
        {\tt cummings@math.huji.ac.il}
        \thanks{Supported by a Postdoctoral Fellowship at the Hebrew University.}
        \\
        Saharon Shelah\\
        Hebrew University of Jerusalem\\
        {\tt shelah@math.huji.ac.il}
        \thanks{Partially supported by the Basic Research Fund of
        the Israel Academy of Science. Paper number 541.}}

\newcommand{\fb}{{\frak b}}
\newcommand{\fd}{{\frak d}}
\newcommand{\fs}{{\frak s}}

\newcommand{\FM}{{\Bbb M}}
\renewcommand{\ll}{{{}^\gl \gl}}

\begin{document}


\baselineskip=16pt
\binoppenalty=10000
\relpenalty=10000

\maketitle

\begin{abstract}
   We prove some consistency results about $\fb(\gl)$ and $\fd(\gl)$,
 which are natural generalisations of the cardinal invariants of the continuum $\frak b$ 
 and $\frak d$. We also define invariants $\fb_{\rm cl}(\gl)$ and 
 $\fd_{\rm cl}(\gl)$, and prove that almost always  $\fb(\gl) = \fb_{\rm cl}(\gl)$
 and  $\fd(\gl) = \fd_{\rm cl}(\gl)$

\end{abstract}

\section{Introduction}

 The cardinal invariants of the continuum have been extensively studied.
 They are cardinals, typically between $\go_1$ and $2^\go$, whose values
 give structural information about ${}^\go \go$. The survey paper
 \cite{JVM} contains a wealth of information about these cardinals.

 In this paper we study some natural generalisations to higher cardinals.
 Specifically, for $\gl$ regular, we define cardinals $\fb(\gl)$ and $\fd(\gl)$ which
 generalise the well-known invariants of the continuum $\fb$ and $\fd$.

 For a fixed value of $\gl$, we will prove that there are some
 simple constraints on the triple of cardinals $(\fb(\gl), \fd(\gl), 2^\gl)$.
 We will also prove that any triple of cardinals obeying these constraints
 can be realised.

 We will then prove that there is essentially no correlation between the
 values of the triple $(\fb(\gl), \fd(\gl), 2^\gl)$ for different values of $\gl$,
 except the obvious one
 that $\gl \longmapsto 2^\gl$ is non-decreasing.
 This generalises Easton's celebrated  theorem (see \cite{Easton}) on the possible
 behaviours of $\gl \longmapsto 2^\gl$;
 since his model was built using Cohen forcing, one can show that in that model $\fb(\gl) = \gl^+$
 and $\fd(\gl) = 2^\gl$ for every $\gl$.

 $\fb(\gl)$ and $\fd(\gl)$ are defined using the co-bounded filter on $\gl$, and
 for $\gl > \go$ we can replace the co-bounded filter by the club filter to
 get invariants $\fb_{\rm cl}(\gl)$ and $\fd_{\rm cl}(\gl)$. We finish the
 paper by proving that these invariants are essentially the same as those
 defined using the co-bounded filter.

 Some investigations have been made into generalising the other 
 cardinal invariants of the continuum, for example in \cite{Z} Zapletal 
 considers $\fs(\gl)$ which is a generalised version of the splitting number
 $\fs$. His work has a different flavour to ours,  since getting $\fs(\gl) > \gl^+$
 needs large cardinals.

 We are indebted to the referee and Jindrich Zapletal for pointing out
a serious  problem with the first version of this paper.

\section{Definitions and elementary facts}

   It will be convenient to define the notions of ``bounding number'' and
 ``dominating number'' in quite a general setting.
  To avoid some trivialities, all partial orderings $\FP$ mentioned in
 this paper (with the exception of the notions of forcing) will
 be assumed to have the property that $\forall p \in \FP \; \exists q
\in \FP \; p <_\FP q$. 
 
\begin{definition} Let $\FP$ be a partial ordering. Then 
\begin{itemize}
\item  $U \subseteq \FP$ is {\it unbounded} if and only if
 $\forall p \in \FP \; \exists q \in U \; q \nleq_\FP p$. 
\item  $D\subseteq \FP$ is {\it dominating} if and only if
 $\forall p \in \FP \; \exists q \in D \; p \le_\FP q$.
\item  $\fb(\FP)$ is the least cardinality of an unbounded subset
 of $\FP$.
\item  $\fd(\FP)$ is the least cardinality of a dominating subset
 of $\FP$.
\end{itemize}
\end{definition}

   The next lemma collects a few elementary facts about the cardinals
 $\fb(\FP)$ and $\fd(\FP)$.

\begin{lemma} 
 Let $\FP$ be a partial ordering, and suppose that $\gb = \fb(\FP)$ and
  $\gd = \fd(\FP)$ are infinite.
Then
\[
    \gb = \cf(\gb) \le \cf(\gd) \le \gd \le \card{\FP}.
\]
\end{lemma}

\begin{proof} To show that $\gb$ is regular, suppose for a contradiction
 that $\cf(\gb) < \gb$. Let $B$ be an unbounded family of cardinality $\gb$,
 and write $B = \bigcup_{\ga < \cf(\gb)} B_\ga$ with $\card{B_\ga} < \gb$.
 For each $\ga$ find $p_\ga$ such that $\forall p \in B_\ga \; p \le p_\ga$,
 then find $q$ such that $\forall \ga < \cf(\gb) \; p_\ga \le q$. Then 
 $\forall p \in B \; p \le q$, contradicting the assumption that $B$ was
 unbounded.

  Similarly, suppose that $\cf(\gd) < \gb$. Let $D$ be dominating with
 cardinality $\gd$ and write $D = \bigcup_{\ga < \cf(\gd)} D_\ga$ where
 $\card{D_\ga} < \gd$. For each $\ga$ find $p_\ga$ such that 
 $\forall p \in D_\ga \; p_\ga \nleq p$, and then find $q$ such that
 $\forall \ga < \cf(\gd) \; p_\ga \le q$. Then $\forall p \in D \; q \nleq p$,
 contradicting the assumption that $D$ was dominating.

\end{proof}

  The next result shows that we cannot hope to say much more.

\begin{lemma} \label{bdposet} Let  $\gb$ and $\gd$ be
 infinite cardinals with
 $\gb = \cf(\gb)$ and $\gd^{<\gb} = \gd$. Define  a partial
 ordering $\FP = \FP(\gb, \gd)$
 in the following way; the underlying set is
 $\gb \times [\gd]^{<\gb}$, and $(\gr, x) \le (\gs, y)$
 if and only if $\gr \le \gs$ and $x \subseteq y$.

  Then $\fb(\FP) = \gb$ and $\fd(\FP) = \gd$.
\end{lemma}

\begin{proof} Let $B \subseteq \FP$ be unbounded. If $\card{B} < \gb$
 then we can define
\begin{eqnarray*}
    \gr & = & \sup \setof{\gs}{\exists y \; (\gs, y) \in B} < \gb, \cr
    x   & = & \bigcup \setof{y}{\exists \gs \; (\gs, y) \in B} \in [\gd]^{<\gb}. \cr
\end{eqnarray*}
 But then $(\gr, x)$ is a bound for $B$, so $\card{B} \ge \gb$ and hence $\fb(\FP) \ge \gb$.
 On the other hand the set $\setof{(\ga, \emptyset)}{\ga < \gb}$ is clearly
 unbounded, so that $\fb(\FP) = \gb$.
         
   Let $D \subseteq \FP$ be dominating. If $\card{D} < \gd$ then
 $\bigcup \setof{y}{\exists \gs \; (\gs, y) \in D} \neq \gd$, and this
 is impossible, so that $\fd(\FP) \ge \gd$.
 On the other hand GCH holds and $\cf(\gd) \ge \gb$,
 so that $\card{\FP} = \gb \times \gd^{<\gb} = \gd$.
 Hence $\fd(\FP) = \gd$.

\end{proof}

\begin{definition} Let $\FP$, $\FQ$ be posets and $f: \FP \lra \FQ$ a function.
 $f$ {\it embeds $\FP$ cofinally into $\FQ$} if and only if 
\begin{itemize}
\item   $\forall p, p' \in \FP \; p \le_\FP p' \iff f(p) \le_\FQ f(p')$.
\item   $\forall q \in \FQ \; \exists p \in \FP \; q \le_\FQ f(p)$. That is,
 $\rge(f)$ is dominating.
\end{itemize}
\end{definition}

\begin{lemma} \label{cofinal} If $f:\FP \lra \FQ$ embeds $\FP$ cofinally into $\FQ$
 then $\fb(\FP) = \fb(\FQ)$ and $\fd(\FP) = \fd(\FQ)$.
\end{lemma}

\begin{proof} Easy.
\end{proof}

\begin{lemma} \label{wflemma} Let $\FP$ be any partial ordering. Then there is
 $\FP^* \subseteq \FP$ such that $\FP^*$ is a dominating subset
 of $\FP$ and $\FP^*$ is well-founded.
\end{lemma}

\begin{proof} 
 We enumerate $\FP^*$ recursively. Suppose that we have already
 enumerated elements $\seq{b_\ga:\ga < \gb}$ into $\FP^*$. If
 $\setof{b_\ga}{\ga<\gb}$ is dominating then we stop, otherwise 
 we choose $b_\gb$ so that $b_\gb \nleq b_\ga$ for all $\ga < \gb$.

 Clearly the construction stops and enumerates a dominating subset $\FP^*$
 of $\FP$. To see that $\FP^*$ is well-founded observe that 
  $b_\gb < b_\ga \implies \gb < \ga$.
\end{proof}

  Notice that the identity embeds $\FP^*$ cofinally in $\FP$, so
 $\fb(\FP^*) = \fb(\FP)$ and $\fd(\FP^*) = \fd(\FP)$.

  We also need some information about the preservation of $\fb(\FP)$ and
 $\fd(\FP)$ by forcing.

\begin{lemma} \label{cclemma} Let $\FP$ be a partial ordering with $\fb(\FP) = \gb$,
 $\fd(\FP) = \gd$. 
\begin{itemize} 
\item Let $V[G]$ be a generic extension of $V$ such that every set
 of ordinals of size less than $\gb$ in $V[G]$ is covered by a set
 of size less than $\gb$ in $V$. Then $V[G] \models \fb(\FP) = \gb$.
\item Let $V[G]$ be a generic extension of $V$ such that every set
 of ordinals of size less than $\gd$ in $V[G]$ is covered by a set
 of size less than $\gd$ in $V$. Then $V[G] \models \fd(\FP) = \gd$.
\end{itemize}
\end{lemma}

\begin{proof} We do the first part, the second is very similar. 
 The hypothesis implies that $\gb$ is a cardinal in $V[G]$, and
 since ``$B$ is unbounded'' is upwards absolute from $V$ to $V[G]$
 it is clear that $V[G] \models \fb(\FP) \le \gb$. Suppose for
 a contradiction that we
 have $C$ in $V[G]$ unbounded with $V[G] \models \card{C} < \gb$.
 By our hypothesis there is $D \in V$ such that $C \subseteq D$
 and $V \models \card{D} < \gb$, but now $D$ is unbounded contradicting
 the definition of $\gb$.
\end{proof}

   With these preliminaries out of the way, we can define the cardinals
 which will concern us in this paper.

\begin{definition} Let $\gl$ be a regular cardinal. 
\begin{enumerate}

\item If $f, g \in \ll$ then $f <^* g$ iff 
 $\exists \ga < \gl \; \forall \gb > \ga \;  f(\gb) < g(\gb)$.

\item $\fb(\gl) =_{def} \fb( (\ll, <^*) )$.

\item $\fd(\gl) =_{def} \fd( (\ll, <^*) )$.


\end{enumerate}
\end{definition}

   These are defined by analogy with some ``cardinal invariants of the
 continuum'' (for a reference on cardinal invariants see \cite{JVM})
 known as $\fb$ and $\fd$. In our notation $\fb = \fb(\omega)$ and
 $\fd = \fd(\omega)$.

\begin{lemma} \label{constraintlemma} 
 If $\gl$ is regular then
\begin{itemize}
\item $\gl^+ \le \fb(\gl)$.
\item $\fb(\gl) = \cf(\fb(\gl))$.
\item $\fb(\gl) \le \cf(\fd(\gl))$.
\item $\fd(\gl) \le 2^\gl$.
\item $\cf(2^\gl) > \gl$.
\end{itemize}
\end{lemma}

\begin{proof} The first claim follows from the following trivial fact.

\begin{fact} Let $\setof{f_\ga}{\ga < \gl} \subseteq \ll$. Then there
 is a function $f \in \ll$ such that $\forall \ga < \gl \; f_\ga <^* f$.
\end{fact}

\begin{proof} Define $f(\gb) = \sup \setof{f_\gg(\gb)+1}{\gg < \gb}$.
 Then $f(\gb) > f_\gg(\gb)$ for $\gg < \gb < \gl$.
\end{proof}

   The next three claims follow easily from our general results on
 $\fb(\FP)$ and $\fd(\FP)$, and the last is just K\"onig's well-known theorem
 on cardinal exponentiation.
   
\end{proof}

  We will prove that these are essentially the only restrictions
 provable in ZFC. One could view this as a refinement of Easton's
 classical result (see \cite{Easton}) on $\gl \longmapsto 2^\gl$.

\section{Hechler forcing}

\label{Hechlerforcing}

   In this section we show how to force that certain posets
 can be cofinally embedded in $(\gl^\gl, <^*)$.  This is
 a  straightforward generalisation of Hechler's work
 in \cite{Hechler}, where he treats the case $\gl = \go$.

 We start with a brief review of our forcing notation.
 $p \le q$ means that $p$ is stronger than $q$, a
 $\gk$-closed forcing notion is one in which every decreasing
 chain of length less than $\gk$ has a lower bound, and
 a $\gk$-dense forcing notion is one in which every sequence
 of dense open sets of length less than $\gk$ has non-empty
 intersection.

 If $x$ is an ordered pair then $x_0$ will denote the first
 component of $x$ and $x_1$ the second component.

\begin{definition} 
 Let $\gl$ be regular. $\FD(\gl)$ is the notion of forcing whose
 conditions are pairs $(s, F)$ with $s \in {}^{<\gl} \gl$ and
 $F \in \ll$, ordered as follows;
  $(s, F) \le (t, F')$ if and only if
\begin{enumerate}
\item $\dom(t) \le \dom(s)$ and $t = s \restriction \dom(t)$.
\item $s(\ga) \ge F'(\ga)$ for $\dom(t) \le \ga < \dom(s)$.
\item $F(\ga) \ge F'(\ga)$ for all $\ga$.
\end{enumerate}
\end{definition}

  We will think of a generic filter $G$ as adding a function $f_G: \gl \lra \gl$
 given by $f_G = \bigcup \setof{s}{\exists F \; (s, F) \in G}$. It is easy to see
 that 
\[
    G = \setof{(t, F)}{t = f_G \restriction \dom(t), \dom(t) \le \ga \implies F(\ga) \le f_G(\ga)},
\]
 so that $V[f_G] = V[G]$ and we can talk about functions from $\gl$ to $\gl$ being
 $\FD(\gl)$-generic.

\begin{lemma} 
 Let $\gl^{<\gl} = \gl$, and set $\FP = \FD(\gl)$. Then
\begin{enumerate}
\item $\FP$ is $\gl$-closed.
\item $\FP$ is  $\gl^+$-c.c.
\item If $g: \gl \lra \gl$ is $\FP$-generic over $V$ then 
 $\forall f \in \ll \cap V \; f <^* g$.
\end{enumerate}
\end{lemma}

\begin{proof}

\begin{enumerate}
\item Let $\gg < \gl$ and suppose that $\seq{(t_\ga, F_\ga): \ga < \gg}$ is
 a descending $\gg$-sequence of conditions from $\FP$. Defining
 $t = \bigcup \setof{t_\ga}{\ga < \gg}$ and $F: \gb \longmapsto \sup \setof{F_\ga(\gb)}{\ga < \gg}$,
 it is easy to see that $(t, F)$ is a lower bound for the sequence. 

\item Observe that if $(s, F)$ and $(s, F')$ are two conditions with the same
 first component then they are compatible, because if $H: \gb \longmapsto F(\gb) \cup F'(\gb)$
 the condition $(s, H)$ is a common lower bound. There are only $\gl^{<\gl} = \gl$ possible
 first components, so that $\FP$ clearly has the $\gl^+$-c.c.

\item Let $f \in \ll \cap V$, and let $(t, F)$ be an arbitrary
 condition. Let us define 
 $F': \gb \longmapsto (F(\gb) \cup f(\gb)) + 1$, then $(t, F')$
 refines $(t, F)$
 and forces that $f(\ga) < f_G(\ga)$ for all $\ga \ge \dom(t)$.

\end{enumerate}
\end{proof}

   If $\gm = \cf(\gm) > \gl = \gl^{<\gl}$ then it is straightforward to iterate
 $\FD(\gl)$ with $<\gl$-support for $\gm$ steps and get a model where
 $\fb(\gl) = \fd(\gl) = \gm$. Getting a model where $\fb(\gl) < \fd(\gl)$  is
 a little harder, but we can do it by a ``nonlinear iteration'' which will embed a
 well chosen poset cofinally into $(\ll, <^*)$.

\begin{theorem} Let $\gl = \gl^{<\gl}$, and suppose that $\FQ$ is any 
 well-founded poset with $\fb(\FQ) \ge \gl^+$. Then there is a forcing 
 $\FD(\gl, \FQ)$ such that
\begin{enumerate} 
\item  $\FD(\gl, \FQ)$ is $\gl$-closed and $\gl^+$-c.c.
\item  $V^{\FD(\gl, \FQ)} \models \hbox{$\FQ$ can be cofinally embedded into $(\ll, <^*)$}$
\item  If $V \models \fb(\FQ) = \gb$ then $V^{\FD(\gl, \FQ)} \models \fb(\gl) = \gb$.
\item  If $V \models \fd(\FQ) = \gd$ then $V^{\FD(\gl, \FQ)} \models \fd(\gl) = \gd$.
\end{enumerate}
\end{theorem}

\begin{proof} We will define the conditions and ordering for  $\FD(\gl, \FQ)$ by
 induction on $\FQ$. The idea is to iterate $\FD(\gl)$ ``along $\FQ$'' so as to
 get a cofinal embedding of $\FQ$ into $\ll$. It will be convenient to 
 define a new poset $\FQ^+$ which consists of $\FQ$ together with a new
 element $top$ which is greater than all the elements of $\FQ$.

 We will define for each $a \in \FQ^+$ a notion of forcing $\FP_{a}$.
 If $a \in \FQ^+$ then we will denote $\setof{c \in \FQ}{c < a}$
 by $\FQ / a$. It will follow from the definition that if $c < a$
 then $\FP_c$ is a complete subordering of $\FP_a$, and that the
 map $p \in \FP_a \longmapsto p \restriction \FQ/c$ is a projection
 from $\FP_a$ to $\FP_c$.

 Suppose that for all $b <_\FQ a$ we have already defined $\FP_{b}$.
 
\begin{enumerate}
\item $p$ is a condition in $\FP_{a}$ if and only if
\begin{enumerate}
\item $p$ is a function, $\dom(p) \subseteq \FQ / a$
 and $\card{\dom(p)} < \gl$. 
\item For all $b \in \dom(p)$, $p(b) = (t, \dot F)$ where $t \in {}^{<\gl} \gl$ and
 $\dot F$ is a $\FP_{b}$-name for a member of $\ll$.
\end{enumerate}
\item If $p, q \in \FP_{a}$ then $p \le q$ if and only if
\begin{enumerate}
\item $\dom(q) \subseteq \dom(p)$.
\item For all $b \in \dom(q)$, if $p(b) = (s, \dot H)$ and $q(b) = (t,\dot I)$ then
\begin{enumerate}
\item   $t = s \restriction \dom(t)$. 
\item   $p \restriction (\FQ / b) \forces_{\FP_{b}} \dom(t) \le \ga < \dom(s) \implies
  s(\ga) > \dot I(\ga)$.
\item   $p \restriction (\FQ / b) \forces_{\FP_{b}} \forall \ga \; \dot H(\ga) \ge \dot I(\ga)$.
\end{enumerate}
\end{enumerate}
\end{enumerate}

  We define $\FD(\gl, \FQ) = \FP_{top}$, and verify that this forcing does what we claimed.
 The verification is broken up into a series of claims.

\begin{claim} $\FD(\gl, \FQ)$ is $\gl$-closed.
\end{claim}

\begin{proof} Let $\gg < \gl$ and let $\seq{p_\ga : \ga < \gg}$ be a descending
 $\gg$-sequence of conditions. We will define a new condition $p$ with
 $\dom(p) = \bigcup_\ga \dom(p_\ga)$. For each
 $b \in \dom(p_\ga)$ let $p_\ga(b) = (t_\ga(b), \dot F_\ga(b))$.

 Let $p(b) = (t, \dot F)$ where
  $t = \bigcup \setof{t_\ga(b)}{b \in \dom(p_\ga)}$ and $\dot F(b)$ is a $\FP_{b}$-name
 for the pointwise supremum of $\setof{\dot F_\ga(b)}{b \in \dom(p_\ga)}$.
 Then it it is easy to check that $p$ is a condition and is a lower bound for
 $\seq{p_\ga: \ga < \gg}$

\end{proof}

\begin{claim} $\FD(\gl, \FQ)$ is $\gl^+$-c.c.
\end{claim}

\begin{proof} Let $\seq{p_\ga : \ga < \gl^+}$ be a family of conditions. Since
 $\gl = \gl^{<\gl}$ we may assume that the domains form a $\gD$-system with
 root $r$. We may also assume that for $b \in r$, $p_\ga(b) = (t_b, F_\ga(b))$
 where $t_b$ is independent of $\ga$. It is now easy to see that any two
 conditions in the family are compatible.
\end{proof}

\begin{claim} If $c < a$ then 
 $\FP_c$ is a complete subordering of $\FP_a$, and  the
 map $p  \longmapsto p \restriction \FQ/c$ is a projection
 from $\FP_a$ to $\FP_c$.
\end{claim}

\begin{proof}
 This is routine.
\end{proof}

   If $G$ is $\FD(\gl, \FQ)$-generic, then for each $a \in \FQ$ we can define 
 $f^a_G \in \ll \cap V[G]$ by
$f^a_G = \bigcup \setof{t(a)_0}{t \in G}$. It is these functions that will
 give us a cofinal embedding of $\FQ$ into $(\ll \cap V[G], <^*)$,
 via the map $a \longmapsto f^a_G$.

\begin{claim} If $a <_\FQ b$ then $f^a_G <^* f^b_G$.
\end{claim}

\begin{proof} Let $p$ be a condition and let $\dot F$ be the canonical
 $\FP_{b}$-name for $f^a_G$. Refine $p$ to $q$ in the following way;
 $q(c) = p(c)$ for $c \neq b$, and if $p(b) = (t, \dot H)$ then
 $q(b) = (t, \dot I)$ where $\dot I$ names the pointwise maximum
 of $\dot F$ and $\dot H$.

              Then $q$ forces that $f^b_G(\ga)$ is greater than
 $f^a_G(\ga)$ for $\ga \ge \dom(t)$.

\end{proof}

 Notice that by the same proof  $f^b_G$ dominates every function
 in $V^{\FP_{a}}$.

\begin{claim} If $a \nless_\FQ b$ then $f^a_G \nless^* f^b_G$.
\end{claim}

\begin{proof}
   If $b <_\FQ a$ then we showed in the last claim that $f^b_G <^* f^a_G$,
 so we may assume without loss of generality that $b \nless_\FQ a$.

 Let $p$ be a condition and let $\ga < \gl$. Choose
 $\gb$ large enough that  $\{ \dom(p(a)_0), \dom(p(b)_0), \ga \} \subseteq \gb$.
 Let $p(b) = (t, \dot F)$, and find $q \in \FP_{b}$ such that
 $q \le p \restriction (\FQ / b)$ and $q$ decides $\dot F \restriction (\gb + 1)$.

 Let $p_1$ be the condition such that $p_1(c) = p(c)$ if $c \nless_\FQ b$ and
 $p_1(c) = q(c)$ if $c <_\FQ b$. Then $p_1$ refines $p$ and $p_1(a) = p(a)$, $p_1(b) = p(b)$.

 Let $p(a) = (s, \dot H)$.
 Find $r \in \FP_{a}$ such that $r \le  p_1 \restriction (\FQ / a)$
 and $r$ decides $\dot H \restriction (\gb + 1)$.

   Let $p_2$ be the condition such that $p_2(c) = p_1(c)$ if $c \nless_\FQ a$ and
 $p_2(c) = r(c)$ if $c <_\FQ a$. Then $p_2$ refines $p_1$ and $p_2(a) = p(a)$, $p_2(b) = p(b)$.

 Now it is easy to extend $p_2$ to a condition  which forces $f^a_G(\gb) > f^b_G(\gb)$.

\end{proof}

\begin{claim} The map $a \longmapsto f^a_G$ embeds $\FQ$ cofinally into
 $(\ll, <^*)$ in the generic extension by $\FD(\gl, \FQ)$. 
\end{claim}

\begin{proof} We have already checked that the map is order-preserving. It remains
 to be seen that its range is dominating.

 Let $G$ be $\FD(\gl, \FQ)$-generic and let $f \in \ll \cap V[G]$.
 Then $f = (\dot f)^G$ for some canonical name $\dot f$, and by the $\gl^+$-c.c.~we
 may assume that there is $X \subseteq \FQ$ such that $\card{X} = \gl$ and $\dot f$
 only involves conditions $p$ with $\dom(p) \subseteq X$. Now $\fb(\FQ) \ge \gl^+$
 so that we can find $a \in \FQ$ with $X \subseteq \FQ / a$. 

  This implies that $\dot f$ is a $\FP_{a}$-name for a function in $\ll$,
 so that $f <^* f^a_G$ and we are done.
\end{proof}

\begin{claim} If $V \models \fb(\FQ) = \gb$ then $V^{\FD(\gl, \FQ)} \models \fb(\gl) = \gb$.
\end{claim}

\begin{proof} Let $G$  be $\FD(\gl, \FQ)$-generic. By lemma \ref{cofinal} it will suffice  to
 show that $V[G] \models  \fb(\FQ) = \gb$. This follows from lemma \ref{cclemma},
 the fact that $\FD(\gl, \FQ)$ is $\gl^+$-c.c.~and the assumption that
 $\fb(\FQ) \ge \gl^+$.
\end{proof}

\begin{claim} If $V \models \fd(\FQ) = \gd$ then $V^{\FD(\gl, \FQ)} \models \fd(\gl) = \gd$.
\end{claim}

\begin{proof} Exactly like the last claim.
\end{proof}
  This finishes the proof of Theorem 1.
\end{proof}

\section{Controlling the invariants at a fixed cardinal}

\label{module}

   In this section we show how to force that the triple $(\fb(\gl), \fd(\gl), 2^\gl)$
 can be anything ``reasonable'' for a fixed value of $\gl$.

\begin{theorem} Let $\gl = \gl^{<\gl}$ and let GCH hold at all cardinals $\gr \ge \gl$. Let
 $\gb, \gd, \gm$ be cardinals such that 
 $\gl^+ \le \gb = \cf(\gb) \le \cf(\gd)$, $\gd \le \gm$ and $\cf(\gm) > \gl$.

Then there is a forcing $\FM(\gl, \gb, \gd, \gm)$ such that in the generic
 extension $\fb(\gl) = \gb$, $\fd(\gl) = \gd$ and $2^\gl = \gm$.
\end{theorem}

\begin{proof} In $V$ define  $\FQ = \FP(\gb, \gd)$, as in lemma \ref{bdposet}.
 We know that $V \models \fb(\FQ) = \gb$ and $V \models \fd(\FQ) = \gd$.
 Fix $\FQ^*$ a cofinal wellfounded subset of $\FQ$, and then define a new
 well-founded poset $\FR$ as follows.

\begin{definition} The elements of $\FR$ are pairs $(p, i)$ where 
 either $i=0$ and $p \in \gm$ or $i=1$ and $p \in \FQ^*$. $(p, i) \le (q, j)$
 iff  $i=j=0$ and $p \le q$ in $\gm$, $i=j=1$ and $p \le q$ in $\FQ^*$, or
 $i=0$ and $j=1$.
\end{definition}

  Now we set  $\FM(\gl, \gb, \gd, \gm) = \FD(\gl, \FR)$. It is routine to use the closure
 and chain condition to argue that $\FM$ makes $2^\gl = \gm$. Since $\FR$ contains a
 cofinal copy of $\FQ^*$, it is also easy to see that $\FM$ forces $\fb(\gl)=\gb$ and
 $\fd(\gl)=\gd$. 

\end{proof}

\section{A first attempt at the main theorem}

\label{try1}

 We now aim to put together the basic modules as described
 in the previous section, so as to control the function
 $\gl \longmapsto (\fb(\gl), \fd(\gl), 2^\gl)$
 for all regular $\gl$. A naive first attempt would
 be to imitate Easton's construction from \cite{Easton};
 this almost works, and will lead us towards the right
 construction.

 Let us briefly recall the statement and proof of Easton's theorem
 on the behaviour of $\gl \longmapsto 2^\gl$. 

\begin{lemma}[Easton's lemma] If $\FP$ is $\gk$-c.c.~and $\FQ$
 is $\gk$-closed then $\FP$ is $\gk$-c.c.~in $V^\FQ$ and
 $\FQ$ is $\gk$-dense in $V^\FP$. In particular
 ${}^{<\gk} ON \cap V^{\FP \times \FQ} = {}^{<\gk} ON \cap V^\FP$.
\end{lemma}

\begin{theorem}[Easton's theorem] \label{ET} Let  $F:REG \lra CARD$
  be  a class function such that $\cf(F(\gl)) > \gl$ and
 $\gl < \gm \implies F(\gl) \le F(\gm)$. Let GCH hold. Then there is a class
 forcing $\FP$ which preserves cardinals and cofinalities, such that in
 the extension $2^\gl = F(\gl)$ for all regular $\gl$.
\end{theorem}

\begin{proof}[Sketch]
 The ``basic module'' is $\FP(\gl) = Add(\gl, F(\gl))$.
 $\FP$ is the ``Easton product'' of the $\FP(\gl)$, to be more
 precise $p \in \FP$ iff
\begin{enumerate}
\item $p$ is a function with $\dom(p) \subseteq REG$ and $p(\gb) \in \FP(\gb)$
  for all $\gb \in \dom(p)$.
\item For all inaccessible $\gg$, $\dom(p) \cap \gg$ is bounded in $\gg$.
\end{enumerate}
 $\FP$ is ordered by pointwise refinement. There are certain complications
 arising from the fact that we are doing class forcing; we ignore them in this
 sketch.

 If $\gb$ is regular then we may factor $\FP$ as $\FP_{<\gb} \times \FP(\gb) \times \FP_{>\gb}$
 in the obvious way. $\FP_{\ge \gb}$ is always $\gb$-closed.

 It follows from GCH and the $\Delta$-system lemma that if $\gg$ is Mahlo or the
 successor of a regular cardinal then $\FP_{<\gg}$ is $\gg$-c.c. On the other hand, if
 $\gg$ is a non-Mahlo inaccessible or the successor of a singular cardinal, then $\FP_{<\gg}$ is
 in general only $\gg^+$-c.c. 

 In particular for $\gg$ regular $\FP_{\le \gg} = \FP_{<\gg^+}$ is always $\gg^+$-c.c.~so
 that by Easton's lemma ${}^\gg ON \cap V^\FP = {}^\gg ON \cap V^{\FP_{\le \gg}}$. This
 implies that in the end we have only added $F(\gg)$ many subsets of $\gg$.

 It remains to be seen that cardinals and cofinalities are preserved. It will suffice
 to show that regular cardinals remain regular. If $\gg$ is Mahlo or the
 successor of a regular cardinal, then Easton's lemma implies that
 ${}^{<\gg} ON \cap V^\FP = {}^{<\gg} ON \cap V^{\FP_{<\gg}}$, and since
 $\gg$ is regular in $V^{\FP_{<\gg}}$ (by $\gg$-c.c.)~$\gg$ is clearly regular
 in $V^\FP$. 

  Now suppose that $\gg = \gm^+$ for $\gm$ singular. If $\gg$ becomes singular in
 $V^\FP$ let its new cofinality be $\gb$, where we see that $\gb < \gm$ and
 $\gb$ is regular in $V$. ${}^\gb ON \cap V^\FP = {}^\gb ON \cap V^{\FP_{\le \gb}}$,
 so that $\gg$ will have cofinality $\gb$ in $V^{\FP_{\le \gb}}$. This is absurd
 as $\FP_{\le \gb}$ is $\gb^+$-c.c.~and $\gb^+ < \gm < \gg$.
 A very similar argument will work in case $\gg$ is a non-Mahlo inaccessible.

\end{proof}

  Suppose that we replace $Add(\gl, F(\gl))$ by 
 $\FP(\gl) = \FM(\gl, \gb(\gl), \gd(\gl), \gm(\gl) )$, where 
 $\gl \longmapsto (\gb(\gl), \gd(\gl), \gm(\gl) )$ is a function obeying the
 constraints given by Lemma \ref{constraintlemma}. Let $\FP$ be the Easton
 product of the $\FP(\gl)$. Then exactly as in the proof of Easton's theorem
 it will follow that $\FP$ preserves cardinals and cofinalities, and that
 $2^\gl = \gm(\gl)$ in $V^\FP$.

\begin{lemma} \label{easycase}
 If $\gl$ is inaccessible or the successor of a regular cardinal then
  $\fb(\gl) = \gb(\gl)$, $\fd(\gl) = \gd(\gl)$ and $2^\gl = \gm(\gl)$ in $V^\FP$.
\end{lemma}

\begin{proof} For any $\gl$, ${}^\gl \gl \cap V^\FP = {}^\gl \gl \cap V^{\FP_{\le \gl}}$.
 $\fb(\gl)$ and $\fd(\gl)$ have the right values in $V^{\FP(\gl)}$ by design, and these
 values are not changed by $\gl$-c.c.~forcing.
 So assuming $\FP_{<\gl}$ is $\gl$-c.c.~those invariants have the right values
 in $V^{\FP_{\le \gl}}$, and hence in $V^\FP$.
\end{proof}

   We need some way of coping with the successors of singular cardinals and the
 non-Mahlo inaccessibles. Zapletal pointed out that at the first inaccessible 
 in an Easton iteration we are certain to add many Cohen subsets, so that there
 really is a need to modify the construction. 

 \section{Tail forcing}

  Easton's forcing to control $\gl \longmapsto 2^\gl$ can be seen as
 a kind of iterated forcing in which we choose each iterand from the
 ground model, or equivalently as a kind of product forcing.
 Silver's ``Reverse Easton forcing'' is an iteration
 in which the iterand at $\gl$ is defined in $V^{\FP_\gl}$.
 The ``tail forcing'' which we describe here is a sort of hybrid.

 We follow the conventions of Baumgartner's  paper \cite{JB}
 in our treatment of iterated forcing, except that when have
 $\dot \FQ \in V^\FP$ and form $\FP * \dot \FQ$ we reserve the
 right not to take all $\FP$-names for members of $\FQ$
 (as long as we take enough names that the set of their denotations
 is forced to be dense).
 For example
 if $\FQ \in V$ we will only take names $\hat q$ for $q \in \FQ$,
 so $\FP * \hat\FQ$ will just be $\FP \times \FQ$.

 We will describe a kind of iteration which we call ``Easton tail iteration''
 in which at successor stages we choose iterands from $V$, but at limit
 stage $\gl$ we choose $\dot \FQ_\gl$ in a different way;  possibly $\FQ_\gl \notin V$,
 but we will arrange things so that the generic $G_\gl$ factors at many places
 below $\gl$ and any final segment of $G_\gl$ essentially determines
 $\FQ_\gl$.
   This idea comes from Magidor and Shelah's paper \cite{MgSh}.
 
 We assume for simplicity that in the ground model all limit cardinals are singular
 or inaccessible. In the application that we intend this is no restriction, as the
 ground model will obey GCH.

\begin{definition}
\label{taildef}
 A forcing iteration $\FP_\gg$ with iterands $\seq{\dot \FQ_\gb: \gb +1 < \gg}$
is an {\em Easton tail iteration\/} iff
\begin{enumerate}
\item The iteration has Easton support, that is to say a direct limit is taken at inaccessible
 limit stages and an inverse limit elsewhere.
\item $\dot \FQ_\gb = 0$ unless $\gb$ is a regular cardinal.
\item If $\gb$ is the successor of a regular cardinal then $\FQ_\gb \in V$.
\item For all regular $\gb$, $\FP_{\gb+1}$ is $\gb^+$-c.c.
\item For $\gl$ a limit cardinal $\FP_\gl$ is $\gl^{++}$-c.c. if $\gl$ is singular,
 and $\gl^+$-c.c. if $\gl$ is inaccessible.
\item For all regular $\ga$ with $\ga + 1 < \gg$ there exists an iteration $\FP^\ga_\gg$ dense
 in $\FP_\gg$ such that $\FP^\ga_{\ga+1} = \FP_{\ga+1}$, and
 for $\gb$ with  $\ga+1 < \gb \le \gg$
\begin{enumerate}
\item  $\FP^\ga_\gb$ factors as $\FP_{\ga+1} \times \FP^\ga \restriction (\ga+1, \gb)$.
\item  If $\gb$ is inaccessible or the successor of a singular, 
 and $p \in \FP^\ga_\gg$, then $p(\gb)$ is a name depending only
 on $\FP^\ga \restriction (\ga+1, \gb)$.
\item  $\FP^\ga \restriction (\ga+1, \gb)$ is $\ga^+$-closed.
\end{enumerate}
\end{enumerate}
\end{definition}

  Clause 6 is of course the interesting one. It holds in a trivial
 way if $\FP_\gg$ is just a product with Easton supports. Clauses 6a
 and 6b should really be read together, as the factorisation in 6a
 only makes sense because 6b already applies to  $\bar\gb < \gb$,
 and conversely 6b only makes sense once we have the factorisation 
 from 6a.

 The following result shows that Easton tail iterations do not disturb the
 universe too much.

\begin{lemma}
 Let $\FP_\gg$ be an Easton tail iteration. Then
\begin{enumerate}
\item For all regular $\ga < \gg$,
 ${}^\ga ON \cap V^{\FP_\gg} = {}^\ga ON \cap V^{\FP_{\ga+1}}$.
\item $\FP_\gg$ preserves all cardinals and cofinalities.
\end{enumerate}
\end{lemma}

\begin{proof} Exactly like Theorem \ref{ET}.
\end{proof}

  In the next section we will see how to define a non-trivial Easton
 tail iteration. If $\gg^+$ is the successor of a regular then it
 will suffice to choose $\FQ_{\gg^+} \in V$ as any $\gg^+$-closed
 and $\gg^{++}$-c.c. forcing. The interesting (difficult) stages are
 the ones where we have to cope with the other sorts of regular
 cardinal, here we will have to maintain the hypotheses on the chain
 condition and factorisation properties of the iteration. It turns out that
 slightly different strategies are appropriate for  inaccessibles
 and successors of singulars.

\section{The main theorem}

\begin{theorem}

\label{mainthm}

 Let GCH hold. Let $\gl \longmapsto (\gb(\gl), \gd(\gl), \gm(\gl))$
 be a class function from $REG$ to $CARD^3$, with 
 $\gl^+ \le \gb(\gl) = \cf(\gb(\gl)) \le \cf(\gd(\gl)) \le \gd(\gl) \le \gm(\gl)$
 and $\cf(\gm(\gl)) > \gl$ for all $\gl$.

 Then there exists a class forcing $\FP_\infty$, preserving all cardinals and cofinalities,
 such that in the generic extension $\fb(\gl) = \gb(\gl)$, $\fd(\gl) = \gd(\gl)$ and
 $2^\gl = \gm(\gl)$ for all $\gl$.
\end{theorem}

\begin{proof}
  We will define by induction on $\gg$ a sequence of Easton tail iterations $\FP_\gg$,
 and then let take a direct limit to get a class forcing $\FP_\infty$. The proof
 that $\FP_\infty$ has the desired properties is exactly as in \cite{Easton}, so we will
 concentrate on defining the $\FP_\gg$. As we define the $\FP_\gg$ we will also 
 define dense subsets $\FP^\ga_\gg$ intended to witness clause 6 in the definition
 of an Easton tail iteration.

 Much of the combinatorics in this section is very similar to that in Section
 \ref{Hechlerforcing}. Accordingly we have only sketched the proofs of some
 of the technical assertions about closure and chain conditions.

  The easiest case to cope with is that where we are looking at the successor of a regular
 cardinal. So let $\gg$ be regular and assume that we have defined $\FP_{\gg^+}$ (which
 is equivalent to $\FP_{\gg+1}$ since we do trivial forcing at all points between $\gg$
 and $\gg^+$), and $\FP^\ga_{\gg^+}$ (which is equivalent to $\FP^\ga_{\gg+1}$) for
 all $\ga \le \gg$.

\begin{definition} $\FQ_{\gg^+} = \FM(\gg^+, \gb(\gg^+), \gd(\gg^+), \gm(\gg^+))$
 as defined in Section \ref{module}. $\FP_{\gg^+ + 1} = \FP_{\gg^+} \times \FQ_{\gg^+}$,
 and  $\FP^\ga_{\gg^+ + 1} = \FP^\ga_{\gg^+} \times \FQ_{\gg^+}$ for $\ga \le \gg$.
\end{definition}

 It is now easy to check that this definition maintains the conditions for being
 an Easton tail iteration. Since $\FP_{\gg^+}$ is $\gg^+$-c.c.~and we know that
 ${}^{\gg^+}\gg^+ \cap V^{\FP_\infty} = {}^{\gg^+}\gg^+ \cap V^{\FP_{\gg^+ + 1}}$,
 we will get the desired behaviour at $\gg^+$ in $V^{\FP_\infty}$.

 Next we consider the case of a cardinal $\gl^+$, where $\gl$ is singular. Suppose
 we have defined $\FP_{\gl^+}$ (that is $\FP_\gl$) and $\FP^\ga_{\gl^+}$ appropriately.
  Let $\FR$
 be a well-founded poset of cardinality $\gm(\gl)$ with $\fb(\FR) = \gb(\gl)$ and
 $\fd(\FR) = \gd(\gl)$, as defined in Section \ref{module}. Let $\FR^+$ be $\FR$ with the
 addition of a maximal element $top$. We will define $\FP_{\gl^+} * \dot \FQ_a$ by induction
 on $a \in \FR^+$, and then set $\FP_{\gl^+ + 1} = \FP_{\gl^+} * \dot \FQ_{top}$.
 In the induction we will maintain the hypothesis that,
for each $\ga < \gl$, $\FP^\ga_{\gl^+} * \dot \FQ_a$ can be factored as
 $\FP_{\ga+1} \times (\FP^\ga \restriction (\ga+1, \gl^+) * \dot \FQ_a)$.
 Let us now fix $a$, and suppose that we have defined $\FP_{\gl^+} * \dot \FQ_b$
 for all $b$ below $a$ in $\FR^+$.

\begin{definition} Let $b < a$, and let $\dot\gt$  be a 
 for a function from $\gl^+$ to $\gl^+$. Then $\dot\gt$ is {\it symmetric} iff for all
 $\ga < \gl$, 
 whenever $G_0 \times G_1$ and $G'_0 \times G_1$ are two generics for
 $\FP_{\ga+1} \times (\FP^\ga \restriction (\ga+1, \gl^+) * \FQ_b)$,
 then $\dot\gt^{G_0 \times G_1} = \dot\gt^{G'_0 \times G_1}$.
\end{definition}

 Of course the (technically illegal) quantification over generic objects in this
 definition can be removed using the truth lemma, to see that the collection of
 symmetric names really is a set in $V$. 

\begin{definition} $(p, q)$ is a condition in $\FP_{\gl^+} * \dot \FQ_a$ iff
\begin{enumerate} 
\item $p \in \FP_{\gl^+}$.
\item $q$ is a function, $\dom(q) \subseteq \FR^+/a$ and $\card{\dom(q)} \le \gl$.
\item For each $b \in \dom(q)$, $q(b)$ is a pair $(s, \dot F)$ where
 $s \in {}^{<\gl^+}\gl^+$ and $\dot F$ is a symmetric $\FP_{\gl^+} * \dot \FQ_b$-name
 for a function from $\gl^+$ to $\gl^+$.
\end{enumerate}
\end{definition}

\begin{definition} Let  $(p, q)$ and $(p', q')$ be conditions  in $\FP_{\gl^+} * \dot \FQ_a$.
$(p', q')$ refines $(p, q)$ iff
\begin{enumerate} 
\item $p'$ refines $p$ in $\FP_{\gl^+}$.
\item $\dom(q) \subseteq \dom(q')$.
\item For each $b \in \dom(q)$, if we let $q(b) = (s, F)$ and $q(b') = (s', F')$,
 then 
\begin{enumerate} 
\item $s'$ extends $s$.
\item If $\ga \in \lh(s') - \lh(s)$ then $(p', q' \restriction (\FR^+/b)) \forces s'(\ga) \ge F(\ga)$.
\item For all $\ga$, $(p', q' \restriction (\FR^+/b)) \forces F'(\ga) \ge F(\ga)$.
\end{enumerate}
\end{enumerate}
\end{definition}

\begin{definition}  We define  $\FP_{\gl^+ +1}$ as  $\FP_{\gl^+} * \dot \FQ_{top}$. If $\ga < \gl$
 then we define $\FP^\ga_{\gl^+ + 1}$ as  $\FP^\ga_{\gl^+} * \dot \FQ_{top}$.
\end{definition}
 
 It is now  routine to check that this definition satisfies the chain condition
 and factorisation demands from  Definition \ref{taildef}.
 The chain condition argument 
 works because $2^\gl = \gl^+$, and any incompatibility in $\FQ_{top}$ is caused
 by a disagreement in the first coordinate at some $b \in \FR$. The factorisation condition
 (clause 6a)
 holds because symmetric names can be computed using any final segment of the
 $\FP_\gl$-generic, and the closure condition (clause 6c) follows from the fact that the canonical
 name for the pointwise sup of a series of functions with symmetric names is itself
 a symmetric name.

 Now we check that we have achieved the desired effect on the values of
 $\fb(\gl^+)$ and $\fd(\gl^+)$.

\begin{lemma} The function added by   $\FP_{\gl^+} * \dot \FQ_{top}$ at $b \in \FR^+$
 eventually dominates all functions  in ${}^{\gl^+}\gl^+ \cap V^{\FP_{\gl^+} * \dot \FQ_b}$.
\end{lemma}

\begin{proof} It suffices to show that if $\dot F$ is a $\FP_{\gl^+} * \dot \FQ_b$-name
 for a function from $\gl^+$ to $\gl^+$ then there is a symmetric name $\dot F'$
 such that $\forces \forall \gb \; F(\gb) \le F'(\gb)$.

  For each regular $\ga < \gl$ we factor $\FP_{\gl^+} * \dot \FQ_b$ as $\FP_{\ga+1}
 \times (\FP^\ga \restriction (\ga+1, \gl^+) * \dot \FQ_b)$.
 In $V^{\FP^\ga \restriction (\ga+1, \gl^+) * \dot \FQ_b }$
 we may treat $\dot F$ as a $\FP_{\ga+1}$-name and define
\[
    G_\ga(\gb) = \sup( \setof{\gg}{\exists p \in \FP_{\ga+1} \; p \forces \dot F(\gb) =\gg}).
\]
 Notice that we may also treat $G_\ga$ as a $\FP_{\gl^+} * \dot \FQ_b$-name, and that
 if $\ga < \bar\ga$ then $\forces G_\ga \le G_{\bar\ga}$. Now let $F'$ be the canonical
 name for a function such that $\forall \gb \; F'(\gb) = \sup_{\ga < \gl} G_\ga(\gb)$,
 then it is easy to see that $F'$ is a symmetric name for a function from $\gl^+$
 to $\gl^+$ and that $\forces \forall \gb \; F(\gb) \le F'(\gb)$.
\end{proof}

\begin{lemma} In $V^{\FP_{\gl^+ + 1}}$ there is a copy of $\FR$ embedded
 cofinally into ${}^{\gl^+} \gl^+$.
\end{lemma}

\begin{proof} Let $\dot f$ name a function from $\gl^+$ to $\gl^+$ in $V^{\FP_{\gl^+ + 1}}$.
 As $\FP_{\gl^+ + 1}$ has the $\gl^{++}$-c.c.~we may assume that $\dot f$ only depends
 on $\gl^+$ many coordinates in $\FR$, and hence (since $\fb(\FR) = \gb(\gl^+) \ge \gl^{++}$)
 that $\dot f$ is a $\FP_{\gl^+} * \dot \FQ_b$-name for some $b \in \FR$. By the preceding
 lemma the function which is added at coordinate $b$ will dominate $\dot f$.
\end{proof}
 
 It remains to be seen what we should do for $\gl$ inaccessible.
 The construction is very similar to that for successors of singulars, with the important
 difference that we need to work with a larger class of names for functions in order
 to guarantee that we dominate everything that we ought to. This in turn leads to a
 slight complication in the definition of $\FP^\ga$.

 Suppose that we have
 defined $\FP_\gl$ and $\FP^\ga_\gl$ appropriately. Let $\FR$ be a poset with the appropriate
 properties ($\card{\FR} = \gm(\gl)$, $\fb(\FR) = \gb(\gl)$, $\fd(\FR) = \gd(\gl)$)
 and let $\FR^+$ be $\FR$ with a maximal element called $top$ adjoined. As before we
 define $\FP_\gl * \dot \FQ_a$ by induction on $a \in \FR^+$. In the induction we will
 maintain the hypothesis that for each $\ga < \gl$ there is a dense subset of $\FP^\ga_\gl * \dot \FQ_a$
 which factorises as $\FP_{\ga+1} \times (\FP^\ga \restriction (\ga+1, \gl) * \dot \FQ^\ga_a)$.
 Let us fix $a$, and suppose that we have defined everything for all $b$ below $a$.

\begin{definition} $(p, (\gm, q))$ is a condition in $\FP_\gl * \dot \FQ_a$
 iff
\begin{enumerate}
\item $p \in \FP_\gl$.
\item $\gm < \gl$, $\gm$ is regular.
\item $q$ is a function, $\dom(q) \subseteq \FR^+/a$ and $\card{\dom(q)} < \gl$.
\item For each $b \in \dom(q)$, $q(b)$ is a pair $(s, \dot F)$ where
 $s \in {}^{<\gl^+}\gl^+$ and $\dot F$ is a  $\FP^\gm \restriction (\gm+1, \gl) * \dot \FQ^\gm_b$-name
 for a function from $\gl^+$ to $\gl^+$.
\end{enumerate}
\end{definition}

\begin{definition} Let  $(p, (\gm, q))$ and $(p', (\gm',q'))$ be conditions  in $\FP_{\gl^+} * \dot \FQ_a$.
$(p', (\gm',q'))$ refines $(p, (\gm,q))$ iff
\begin{enumerate} 
\item $p'$ refines $p$ in $\FP_{\gl^+}$.
\item $\dom(q) \subseteq \dom(q')$.
\item $\gm' \ge \gm$.
\item For each $b \in \dom(q)$, if we let $q(b) = (s, F)$ and $q(b') = (s', F')$,
 then 
\begin{enumerate} 
\item $s'$ extends $s$.
\item If $\ga \in \lh(s') - \lh(s)$ then $(p', q' \restriction (\FR^+/b)) \forces s'(\ga) \ge F(\ga)$.
\item For all $\ga$, $(p', q' \restriction (\FR^+/b)) \forces F'(\ga) \ge F(\ga)$.
\end{enumerate}
\end{enumerate}
\end{definition}

  Notice that since any $\FP_\gl * \dot \FQ_b$-generic induces a 
 $\FP^\gm \restriction (\gm+1, \gl) * \dot \FQ^\gm_b$-generic, there is a natural
 interpretation of any  $\FP^\gm \restriction (\gm+1, \gl) * \dot \FQ^\gm_b$-name as
 a  $\FP_\gl * \dot \FQ_b$-name. We are using this fact implicitly when we define
 the ordering on the conditions. We need to maintain the hypothesis on factorising the
 forcing, so we make the following definition.

\begin{definition} Let $\gm$ be regular with $\gm < \gl$.  
 $\FP_{\gm+1} \times (\FP^\gm \restriction (\gm+1, \gl) * \dot \FQ^\gm_a)$
 is defined as the set of $(p_0, (p_1, (\gn, q)))$ such that
$p_0 \in \FP_{\gm+1}$, $p_1 \in \FP^\gm \restriction (\gm+1, \gl)$
 and $(p_0 \frown p_1, (\gn, q)) \in \FP_\gl * \dot \FQ_a$ with  $\gn \ge \gm$.
\end{definition}

 The key point here is that the factorisation makes sense, because for such a
 condition $q(b)$ depends  only on $\FP^\gm \restriction (\gm+1, \gl) * \dot \FQ^\gm_b$.

\begin{definition} We define $\FP_{\gl+1}$ as $\FP_\gl * \dot \FQ_{top}$,
 and $\FP^\gm_{\gl+1}$ as $\FP_\gl * \dot \FQ^\gm_{top}$, 
\end{definition}

    As in the case of a successor of a singular, it is straightforward to see that we have satisfied
 the chain condition and factorisation conditions. To finish the proof we need to
 check that the forcing at $\gl$ has achieved the right effect, which will be clear
 exactly as in the singular case when we have proved the following lemma.

\begin{lemma} The function added by $\FP_\gl * \dot \FQ_{top}$ at $b \in \FR$ eventually
 dominates all functions in ${}^\gl \gl \cap V^{\FP_\gl * \dot \FQ_b}$. 
\end{lemma}

\begin{proof} We do a density argument. Suppose that $\dot f$ names a function in
 ${}^\gl \gl \cap V^{\FP_\gl * \dot \FQ_b}$, and let $(p, (\gm, q))$ be a condition
 in $\FP_\gl * \dot \FQ_{top}$. 
 We factor $\FP_\gl * \dot \FQ_b$ as 
$\FP_{\gm+1} \times (\FP^\gm \restriction (\gm+1, \gl) * \dot \FQ^\gm_b)$,
 and use the fact that $\FP_{\gm+1}$ has $\gm^+$-c.c.~in
 $V^{\FP^\gm \restriction (\gm+1, \gl) * \dot \FQ^\gm_b}$ to find a 
 $\FP^\gm \restriction (\gm+1, \gl) * \dot \FQ^\gm_b$-name $\dot g$
 such that $\forces \forall \gb \; \dot f(\gb) \le \dot g(\gb)$.
 
  Now we can refine $(p, (\gm, q))$ in the natural way by strengthening the
 second component of $q(b)$ to dominate $\dot g$. This gives a condition which forces that the
 function added at $b$ will eventually dominate $\dot f$.
\end{proof}
  This concludes the proof of Theorem \ref{mainthm}.
\end{proof}

\section{Variations}
\newcommand{\lcl}{{<_{\rm cl}}}

 In this appendix we discuss the invariants that arise if we work
 with the club filter in place of the co-bounded filter.
 It turns out that this does not make too much difference.
All the results here are due to Shelah.

\begin{definition} Let $\gl$ be regular. 
\begin{enumerate}

\item Let $f, g \in \ll$. $f \lcl g$ iff there is
 $C \subseteq \gl$ closed and unbounded in $\gl$
 such that $\ga \in C \implies f(\ga) < g(\ga)$.
\item $\fb_{\rm cl}(\gl) =_{def} \fb( (\ll, \lcl) )$.
\item $\fd_{\rm cl}(\gl) =_{def} \fd( (\ll, \lcl) )$.
\end{enumerate}
\end{definition}

\begin{theorem} \label{dthm} $\fd_{\rm cl}(\gl) \le \fd(\gl) \le
\fd_{\rm cl}(\gl)^\go$. 
\end{theorem}

\begin{proof} If a family of functions is dominating with respect to $<^*$
it  is dominating with respect to $\lcl$, so that $\fd_{\rm cl}(\gl)
\le \fd(\gl)$. 

 For the converse, let us fix $D \subseteq \ll$ such that $D$ is dominating 
 with respect to $\lcl$ and $\card{D} = \fd_{\rm cl}(\gl)$.
 We may assume that every function in $D$ is increasing (replace each
 $f \in D$ by $f^*:\gg \longmapsto \bigcup_{\ga \le \gg} f(\ga)$).

 Let $g_0 \in \ll$. Define by induction $f_n$, $g_n$, and $C_n$ such that
\begin{enumerate}
\item $f_n \in D$.
\item $C_n$ is club in $\gl$, and $\ga \in C_n \implies g_n(\ga) < f_n(\ga)$.
\item $C_{n+1} \subseteq C_n$.
\item $g_n$ is increasing.
\item $g_{n+1}(\gb) > g_n(\gb)$ for all $\gb$. 
\item $g_{n+1}(\gb) > f_n(\min(C_n - (\gb+1)))$ for all $\gb$.
\end{enumerate}

 Now let $\ga = \min( \bigcap_{n < \go} C_n )$.
We will prove that
 $g_0(\gg) \le \bigcup_n f_n(\gg)$
for  $\gg > \ga$.

 Fix some $\gg > \ga$. For each $n$ we know that $C_n \cap \gg \neq \emptyset$,
 so that if we define $\gg_n = \sup( C_n \cap (\gg+1))$ 
 then $\gg_n$ is the largest point of $C_n$ less than or equal to $\gg$.
 Notice that $\min(C_n -(\gg+1)) = \min( C_n -(\gg_n+1))$.

 Since $C_{n+1} \subseteq C_n$, $\gg_{n+1} \le \gg_n$, so that for all
 sufficiently large $n$ (say $n \ge N$) we have $\gg_n = \bar \gg$ for some
 fixed $\bar\gg$.
 We claim that $g_0(\gg) \le f_{N+1}(\gg)$, which we will prove
 by building a chain of inequalities.
 Let us define 
\[
\gd = \min(C_N - (\gg+1)) = \min(C_N - (\bar\gg+1)).
\]
 Then 
\[
   g_0(\gg) \le g_N(\gg) \le g_N(\gd) < f_N(\gd) < g_{N+1}(\bar \gg) < f_{N+1}(\bar\gg) \le f_{N+1}(\gg),
\]
 where the key point is that $\bar\gg \in C_{N+1}$ and hence
$g_{N+1}(\bar \gg) < f_{N+1}(\bar\gg)$. 

  Now it is easy to manufacture a family of size $\fd_{\rm cl}(\gl)^\go$ which
 is dominating with respect to $<^*$, so that $\fd(\gl) \le \fd_{\rm cl}(\gl)^\go$.
\end{proof}

\begin{theorem} $\fb_{\rm cl}(\gl) = \fb(\gl)$. 
\end{theorem}

\begin{proof}
 If a family of functions is unbounded with respect to $\lcl$ it is
 unbounded with respect to $<^*$, so that $\fb(\gl) \le \fb_{\rm cl}(\gl)$.
 
 Suppose for a contradiction that $\fb(\gl) < \fb_{\rm cl}(\gl)$, and fix
 $U_0 \subseteq \ll$ such that $\card{U_0} = \fb(\gl)$ and $U_0$ is
 unbounded with respect to $<^*$. We may assume without loss of
 generality that every function in $U_0$ is increasing.
 We perform an inductive construction in $\go$ steps, whose
 aim is to produce a bound for $U_0$ with respect to $<^*$.

 By assumption $U_0$ is bounded with respect to $\lcl$, so choose $g_0$
 which bounds it modulo the club filter. Choose also club sets
 $\setof{C^0_f}{f \in U_0}$ such that
 $\ga \in C^0_f \implies f(\ga) < g_0(\ga)$.
 For each $f \in U_0$ define another function $f^{[0]}$ by
 $f^{[0]} : \gb \longmapsto f(\min(C^0_f - (\gb+1)))$

 For $n \ge 1$ define $U_n = U_{n-1} \cup  \setof{ f^{[n-1]} }{f \in U_{n-1}}$.
 By induction it will follow that $\card{U_n} = \fb(\gl)$, so that  we may choose $g_n$
 such that $g_n(\gb) > g_{n-1}(\gb)$ for all $\gb$ and $g_n$ bounds $U_n$
 modulo clubs. We choose clubs $\setof{C^n_f}{f \in U_n}$ such that  
\begin{enumerate}

\item  $\ga \in C^n_f \implies f(\ga) < g_n(\ga)$.

\item  If $f \in U_{n-1}$, then $C^n_f \subseteq C^{n-1}_f$.

\item  If $n \ge 2$ and  $f \in U_{n-2}$ then $C^n_f \subseteq C^{n-1}_{f^{[n-2]}}$,
 where this makes sense because in this case $f^{[n-2]} \in U_{n-1}$.

\end{enumerate}
 For each $f \in U_n$ we define 
 $f^{[n]} : \gb \longmapsto f(\min(C^n_f - (\gb+1)))$ to finish
 round $n$ of the inductive construction.

   Now we claim that the pointwise sup of the sequence $\seq{g_n : n<\go}$ is
 an upper bound for $U_0$ with respect to $<^*$. Let us fix $f \in U_0$, and
 then let $\ga = \min(\bigcap C^n_f)$.

 We will now give a very similar argument to that of Theorem \ref{dthm}.
 Fix $\gg > \ga$.
 We  define
 $\gg_n = \sup(C^n_f \cap (\gg+1))$ and observe that $\gg_{n+1} \le \gg_n$,
 so we may find $N$ and $\bar\gg$ such that $n \ge N \implies \gg_n = \bar\gg$.

 Let $\gd = \min(C^N_F - (\gg+1)) = \min(C^N_F - (\bar\gg+1))$.
 We now get a chain of inequalities
\[
  f(\gg) \le f(\gd) = f^{[N]}(\bar\gg) < g_{N+1}(\bar\gg) \le g_{N+1}(\gg).
\]
  This time the key point is that $\bar\gg \in C^{N+2}_f \subseteq C^{N+1}_{f^{[N]}}$,
 so that  $f^{[N]}(\bar\gg) < g_{N+1}(\bar\gg)$.

We have proved that  $\gg > \ga \implies f(\gg) \le \bigcup_n
g_n(\gg)$, so that 
 every function $f \in U_0$ is bounded on a final segment of $\gl$ by
 $\gg \longmapsto \bigcup_n g_n(\gg)$. This contradicts the choice
 of $U_0$ as unbounded with respect to $<^*$, so we are done.

\end{proof}

 It is natural to ask whether the first result can be improved to show that
 $\fb_{\rm cl}(\gl) = \fb(\gl)$. This can be done for $\gl$ sufficiently large,
 at the cost of using a powerful result from Shelah's paper \cite{Sh460}.

\begin{definition} Let $\gth = \cf(\gth) < \gm$. 
\begin{enumerate} 
\item $\powerset^\gth(\gm) = \setof{X \subseteq \gm}{\card{X} = \gth}$.
\item $\gm^{[\gth]}$ is the least cardinality of a family $P \subseteq \powerset^\gth(\gm)$
 such that 
\[
   \forall A \in \powerset^\gth(\gm) \; \exists B \subseteq P \;
 (\card{B} < \gth \wedge A \subseteq \bigcup B).
\]
\end{enumerate}
\end{definition}

  One of the main results of \cite{Sh460} is  that ZFC proves a weak
form of the GCH. 

\begin{theorem} \label{heavyvoodoo} Let $\gm > \beth_\go$. Then
$\gm^{[\gth]} = \gm$  for all sufficiently large $\gth < \beth_\go$.
\end{theorem}

  It is easy to see that if $P \subseteq \powerset^\gth \gm$ is such that 
\[
  \forall A \in \powerset^\gth(\gm) \; \exists B \subseteq P \;
 (\card{B} < \gth \wedge A \subseteq \bigcup B),
\]
  then 
$\forall A \in \powerset^\gth(\gm) \; \exists C \in P \; \card{A \cap C}=\gth$.
 This is all we use in what follows, and in fact we could get away with 
$\forall A \in \powerset^\gth(\gm) \; \exists C \in P \; \card{A \cap
C}=\ha_0$. 

\begin{theorem} \label{bthm2} Let $\gl = \cf(\gl) > \beth_\go$. Then
$\fd(\gl) = \fd_{\rm cl}(\gl)$. 
\end{theorem}

\begin{proof} Let $\gm = \fd_{\rm cl}(\gl)$. Then $\gm > \beth_\go$,
 so that we may apply Theorem \ref{heavyvoodoo} to find  a regular
$\gth < \beth_\go$ such that $\gm^{[\gth]} = \gm$. let us fix $P
\subseteq \powerset^\gth (\gm)$ such that $\card{P} = \gm$ and
$\forall A \in \powerset^\gth(\gm) \; \exists C \in P \; \card{A \cap
C}=\gth$.  

 Now let $D \subseteq \ll$ be such that $\card{D} = \gm$ and $D$ is dominating
 in $(\ll, \lcl)$.
 We may suppose that $D$ consists of increasing functions.
 Enumerate $D$ as $\seq{h_\ga:\ga<\gth}$, and then define
 $h_A: \gg \longmapsto \bigcup_{\ga \in A} h_\ga(\gg)$ for each $A \in P$.
 Since $\gth < \beth_\go \le \gl$, $h_A \in \ll$. We will prove that
 $\setof{h_A}{A \in P}$ is dominating in $(\ll, <^*)$.

 We will do a version of the construction from Theorem \ref{dthm}. 
 Let $g_0 \in \ll$.
 Define by induction $f_\ga$, $g_\ga$, and $C_\ga$
 for $\ga < \gth$, with the following properties.
\begin{enumerate}

\item $f_\ga \in D$.

\item $C_\ga$ is club in $\gl$, and $\gb \in C_\ga \implies g_\ga(\gb)
< f_\ga(\gb)$. 

\item $\ga < \bar\ga \implies C_{\bar\ga} \subseteq C_{\ga}$.

\item $g_\ga$ is increasing.

\item If $\ga < \bar\ga$, then 
 $g_{\bar\ga}(\gb) > g_{\ga}(\gb)$ for all $\gb$.
 
\item  If $\ga < \bar\ga$, then 
 $g_{\bar\ga}(\gb) > f_{\ga}(\min(C_{\ga} - (\gb+1)))$ for all $\gb$.
\end{enumerate}

 This is easy, because $\gth < \gl$. By the choice of $P$ we may
 find a set $A \in P$ such that
$\card{ \setof{h_\gb}{\gb \in A} \cap \setof{f_\ga}{\ga < \gth}} = \gth$.
 Enumerate the first $\go$ many $\ga$ such that $f_\ga \in
\setof{h_\gb}{\gb \in A}$ as $\seq{\ga_n:n<\go}$.

  We may now repeat the proof of Theorem \ref{dthm} with $f_{\ga_n}$,
$g_{\ga_n}$ 
 and $C_{\ga_n}$ in place of $f_n$, $g_n$ and $C_n$. We find that for
all sufficiently  large $\gg$ we have
 $g_0(\gg) \le g_{\ga_0}(\gg) \le \bigcup_n f_{\ga_n}(\gg)$.
 By the definition of $h_A$ and the fact that
$\setof{f_{\ga_n}}{n<\go} \subseteq \setof{h_\gb}{\gb \in A}$, 
 $\bigcup_n f_{\ga_n}(\gg) \le h_A(\gg)$ for all $\gg$, so that $g_0 <^* h_A$.
 This shows $\setof{h_A}{A \in P}$ to be dominating, so we are done.
 
\end{proof}
   
  We do not know whether it can ever be the case that $\fd_{\rm
cl}(\gl) < \fd(\gl)$. This is connected with some open and apparently
difficult questions in pcf theory.

\end{document}